\documentclass[%
 aip,
 cha,
 amsmath,amssymb,
 reprint,%
]{revtex4-2}

\usepackage{graphicx}
\usepackage[usenames,dvipsnames]{color}

\usepackage[utf8]{inputenc}
\usepackage[T1]{fontenc}
\usepackage{mathptmx}

\newcommand{\beq}{\begin{equation}}
\newcommand{\eeq}{\end{equation}}
\newcommand{\bseq}{\begin{subequations}}
\newcommand{\eseq}{\end{subequations}}

\newcommand{\rf}     [1] {~\cite{#1}}
\newcommand{\refref} [1] {Ref.~\cite{#1}}

\newcommand{\reffig} [1] {Fig.~\ref{#1}}

\newcommand{\refFig} [1] {Figure~\ref{#1}}

\newcommand{\refeq}  [1] {Eq.~(\ref{#1})}
\newcommand{\refEq}  [1] {Eq.~(\ref{#1})}
\newcommand{\refsect}[1]{Sec.~\ref{#1}}
\newcommand{\refapp}[1]{Appendix~\ref{#1}}
\newcommand{\reftab} [1] {Tab.~\ref{#1}}

\newcommand{\etc}{{etc.}}       

\newcommand{\ie}{{i.e.}}        
\newcommand{\cf}{{\em cf.\ }}   
\newcommand{\eg}{{e.g.\ }}

\newcommand{\poinc}{\ensuremath{\mathcal{P}}} 
\newcommand{\frm}{\ensuremath{R}} 
\newcommand{\descrn}[1]{\ensuremath{_{#1}}} 
\newcommand{\dr}{\ensuremath{d_\mathrm{r}}} 
\newcommand{\df}{\ensuremath{d_\mathrm{f}}} 
\newcommand{\stitle}{Manifold learning reveals the topology of chaotic flows}

\usepackage{hyperref}
\hypersetup{
    pdfauthor=Siminos,
    pdftitle=\stitle
}
\definecolor{hreflinkcolor}{rgb}{0.13,0.17,0.83}
\hypersetup{colorlinks=true,urlcolor=hreflinkcolor,
		linkcolor=hreflinkcolor,citecolor=hreflinkcolor}

\begin{document}

\title[\stitle]{Manifold learning of Poincar\'e sections reveals the topology of high-dimensional chaotic flows}

\author{Evangelos Siminos}\email{evangelos.siminos@gmail.com}
\affiliation{Department of Physics, University of Gothenburg, SE-41296 G\"{o}teborg, Sweden}

\date{\today}

\begin{abstract}
It is shown that applying manifold learning techniques to Poincar\'e sections of high-dimensional, chaotic dynamical systems can uncover their low-dimensional topological organization.  Manifold learning provides a low-dimensional embedding and intrinsic coordinates for the parametrization of data on the Poincar\'e section, facilitating the construction of return maps with well defined symbolic dynamics. The method is illustrated by numerical examples for the R\"ossler attractor and the Kuramoto-Sivashinsky equation. For the latter we present the reduction of the high-dimensional, continuous-time flow to dynamics on one- and two two-dimensional Poincar\'e sections. We show that in the two-dimensional embedding case the attractor is organized by one-dimensional unstable manifolds of short periodic orbits. In that case, the dynamics can be approximated by a map on a tree which can in turn be reduced to a trimodal map of the unit interval. In order to test the limits of the one-dimensional map approximation we apply classical kneading theory in order to systematically detect all periodic orbits of the system up to any given topological length. 
\end{abstract}

\maketitle

\begin{quotation}
When confronted with chaotic dynamical systems, there are different levels at which we may declare our that our aim to understand them has been fulfilled. A classical approach that goes back to the early days of the field is to reduce the continuous time dynamics to a discrete time map on a Poincar\'e section, a hyperplane that intersects trajectories transversally. That approach, when successful, elucidates the topological organization of the flow and allows one to label different solutions with symbolic sequences and to reduce chaotic dynamics to a ``walk'' in the space of such sequences known as \emph{symbolic dynamics}. This enhances our understanding of how chaotic dynamics are generated for a given flow and allows dynamical systems to be classified in equivalence classes. Moreover, it allows us to  systematically generate and classify the compact solutions that organize chaotic dynamics. However, carrying out such a topological program is far from trivial for high-dimensional dynamical systems, such as fluid flows, even when such systems are known to posses low dimensional attractors. The difficulty lies in the lack of a systematic way to find low dimensional coordinates on the attractor that would lead to useful return maps. Here, it is shown that manifold learning techniques, when applied directly on the Poincar\'e section rather than on the continuous time flow, can provide a straightforward way to construct low-dimensional return maps.
\end{quotation}

\section{\label{sec:intro}Introduction}


\begin{figure*}[ht!]
    \includegraphics[width=\textwidth]{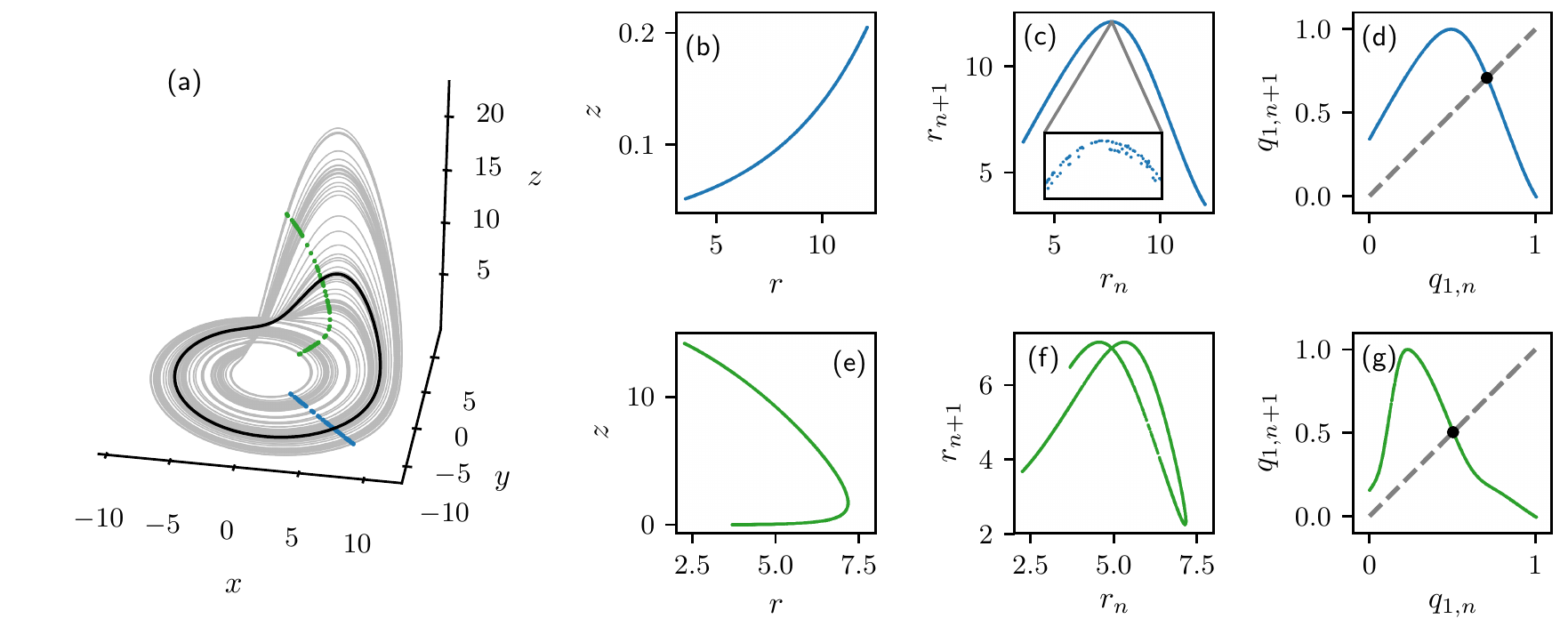}
        \caption{\label{f:rossler_maps} (a) Long trajectory on the attractor of the R\"ossler flow (grey), 
        periodic orbit $\overline{1}$ (black), points of intersection with the Poincar\'e surface of section $\poinc_{1}$ (blue) and $\poinc_{2}$ (green), 
        (b,e) points of intersection with $\poinc_{1}$ (blue) and $\poinc_{2}$ (green) in $(r,z)$ representation, 
        (c,f) corresponding first return maps in the radial variable; the inset in (c) shows an area with dimensions $(\delta r_n,  \delta r_{n+1})=(2.5\times10^{-3}, 1.7\times10^{-6})$,
        (d,g) return maps in intrinsic coordinate $q$ found through application of LLE with parameters $K=160$, $Np=4372$ and $\epsilon=1e-3$, 
        along with the diagonal  $q_{n+1}=s_n$ (dashed grey line) and the periodic orbit $\overline{1}$ (black dot).
        }
\end{figure*}

One of the most successful methods to uncover the topological organization of chaotic attractors is to convert the continuous time flow defined in a state-space of dimension $d$ to a discrete time return map defined on a $d-1$-dimensional hyper-surface known as a Poincar\'e surface of section. Variations of the Poincar\'e section technique have allowed to represent low-dimensional dynamical systems as one- or two-dimensional maps and to apply the well developed methods of topological analysis of discrete maps in order to shed insight in the dynamics of the original dynamical system~\cite{gilmore2003,DasBuch}. In particular, this approach allows the spectrum of 
admissible compact solutions, such as periodic orbits, to be determined and attractors to be classified in topological conjugacy classes. However, this program becomes rather hard  to carry out in high-dimensional dynamical systems, for which there is usually no obvious way to choose and parametrize a Poincar\'e section that would lead to a useful return map. This hinders progress in the study of turbulent fluid flows as high-dimensional dynamical systems for which there has been considerable interest over the last twenty years, see for example~\refref{KawKida01, Visw07b, GHCW07, WFSBC15, Suri20, GraFlo21}. While the state space of fluid flows, described by partial differential equations (PDEs), is formally infinite dimensional, Lyapunov analysis arguments~\rf{TaGiCh11,DCTSCD14} and the theory of inertial manifolds~\rf{Foias1988a,FNSTks88,constantin_integral_1989} suggest that the asymptotic dynamics of dissipative PDEs may take place in a finite-dimensional state-space. 
However, it remains unclear how to construct lower-dimensional representations of fluid flows that could be utilized in dynamical systems studies.

The topic of dimensionality reduction has been pursued in the field of nonlinear dynamics of dissipative PDEs using many different approaches. Linear methods, most notably  proper orthogonal decomposition~\cite{Holmes96} fail to capture the nonlinear structure of the attractor or inertial manifold and are often sub-optimal in terms of number of variables required to faithfully capture the dynamics~\cite{gibsonPhD}. Nonlinear dimensionality reduction, also known as manifold learning~\cite{LeeVer07,Izenman12} (reviewed briefly in \refsect{s:ML}), could be a better alternative since it provides methods to re-embed geometric objects from a high-dimensional to a lower dimensional space, while preserving geometrical or topological structure. While several works, \eg~\cite{Bollt07,Nadler06,Dsilva18}, have applied manifold learning techniques 
to data generated by chaotic dynamical systems, it is still unclear that such methods could help elucidate the topology of high-dimensional dynamical systems, for two reasons. First, such methods need to operate on finely sampled, high-dimensional data in order to correctly reproduce the low-dimensional structure, limiting their scalability towards systems with a high-dimensional state-space. Second, even when a low dimensional embedding of, \eg, a chaotic attractor can be obtained, the underlying topology often remains obscure\rf{Bollt07}.

In this work it is proposed that, for the purpose of reducing a continuous time, high-dimensional dynamical system to low-dimensional return map, nonlinear dimensionality reduction methods should be applied directly on the Poincar\'e section rather than on the continuous time system. While the introduction of a Poincar\'e section is associated with a small reduction in dimensionality (from $d$ to $d-1$ dimensions), it leads to a dramatic reduction in the amount of data required to describe the system, since it eliminates the need to sample trajectories (in time) between successive intersections. As shown in the example of the low-dimensional R\"ossler flow, \refsect{s:method},  manifold learning solves the problem of choosing a suitable surface of section in a natural way, by providing an intrinsic parametrization of (data on) the section. This leads to well defined return maps independently of the choice of surface of section. The method can be seen as a generalization of the idea of using the geodesic distance along ordered periodic points~\cite{Christiansen97} or one-dimensional unstable manifolds of periodic orbits to parametrize the return map~\cite{lanCvit07,SiCvi10}. However, manifold learning is much less restrictive, since it does not require the computation of periodic orbits and their unstable manifolds and can also parameterize multi-dimensional manifolds. This allows us to treat more complicated topologies than hitherto possible.

In order to put this idea to the test for  higher-dimensional dynamical systems, the Kuramoto-Sivashinsky partial differential equation (KSe) is studied for two parameter values corresponding to Poincar\'e section data that can be parametrized as one- or two-dimensional manifolds in \refsect{s:ks_1D} and \refsect{s:ks_2D}, respectively. Surprisingly, we find that even in the case of a two-dimensional embedding, a description based on one-dimensional tree maps is possible. Since such tree maps can be somewhat cumbersome to work with,  we show that they can be reduced to ordinary maps of the unit interval. In order to illustrate the utility but also the limits of the one-dimensional map approximation and the corresponding symbolic dynamics we use classical kneading theory to systematically determine KSe periodic orbits in \refsect{s:ks_cycles}. In \refsect{s:ks_UM} we examine the organization of the attractor in terms of the one-dimensional unstable manifolds of the shortest periodic orbits and argue that the presented method could be useful whenever dynamics is organized by low-dimensional unstable manifolds.

\section{\label{s:method}Manifold learning of Poincar\'e sections}

\subsection{\label{s:rossler}Motivation: Return maps of chaotic flows}

To illustrate the problem of Poincar\'e section choice, and present it in the context of the R\"ossler flow~\rf{ross}, along the lines of \refref{DasBuch}. The R\"ossler flow reads
\begin{align}\label{eq:rossler}
    \dot{x} &=-y-z\,,\\
    \dot{y} &=x+a\,y\,\\
    \dot{z} &=b+z\,(x-c)\,.
\end{align}
In the following we choose parameters $a=b=0.2$, $c=5.7$, for which \refeq{eq:rossler} has a chaotic attractor with Kaplan-York estimate~\cite{KapYor79a} of the fractal dimension $d_f\simeq2.03$, see \refref{KuMoVa14}. \refFig{f:rossler_maps}(a) shows a state-space portrait obtained by integration of a random initial condition which converges, after a short transient, to the strange attractor. 

For a continuous-time dynamical system with $d$-dimensional state-space, 
a Poincar\'e surface of section can be defined as a $(d-1)$-dimensional hypersurface 
which intersects the flow transversally, supplemented by a condition for the direction of the crossing. 
For the R\"ossler flow a natural choice is a two-dimensional plane orthogonal to the $(x,y)$-plane. 
Introducing cylindrical coordinates $(r,\theta,z)$, where $r^2=x^2+y^2$ and $\theta=\arctan(y/x)$, 
we illustrate two different Poincar\'e section choices $\poinc_{1}$, $\poinc_{2}$ 
defined by $\theta=-\pi/4$ and $\theta=\pi/4$, respectively, in \reffig{f:rossler_maps}(a). 
The  corresponding points of intersection are shown in the two-dimensional plane $(r,z)$ 
in \reffig{f:rossler_maps}(b,e) and can be used to define a \emph{return map} $\frm$ 
of successive intersections of the flow with the Poincar\'e section 
$\frm_{}:\ (r\descrn{n},z\descrn{n})\mapsto(r\descrn{n+1},z\descrn{n+1})$, 
where $r\descrn{n}$ denotes the value of $r$ at the $n$'th intersection, \etc.  
Since the points of intersection lie on an approximately one-dimensional curve 
one may attempt to reduce the dimension of the mapping by using a single coordinate, 
for example $r$, to construct a return map  $\frm_{}:\ r\descrn{n}\mapsto r\descrn{n+1}$. Such return maps $\frm_{1}$, $\frm_{2}$ corresponding to the sections $\poinc_{1}$, $\poinc_{2}$, respectively, are shown in \reffig{f:rossler_maps}(c,f). 

The unimodal (\ie, having a single maximum) return map $\frm_{1}$ is approximately single valued, up to 
the ``transverse structure'' or ``thickness'' revealed in the inset of \reffig{f:rossler_maps}(c). The origin of this structure lies in the non-integer dimension of the attractor. In what follows we ignore this fine structure and consider such return maps to be one-dimensional and smooth. This approximation allows to apply the extensive set of tools developed for one-dimensional maps on the unit interval to return maps obtained from continuous-time flows by the Poincar\'e section method~\cite{gilmore2003,DasBuch}.
By assigning the symbols 0 and 1 to points  on the left and right, respectively, of the critical point $r_c$
(the point that corresponds to the maximum of $\frm_1$), trajectories can be encoded as symbolic sequences or \emph{itineraries}. 
We will forgo presenting the topological program in terms of the exhaustively studied R\"ossler system 
and come back to this point in the discussion of Kuramoto-Sivashinsky equation in~\refsect{s:ks_cycles}, where we will also see that ignoring the transverse structure imposes limitations on the validity of this symbolic encoding for longer orbits.

By contrast to $\frm_{1}$, the return map $\frm_{2}$ obtained using the section $\poinc_2$, \reffig{f:rossler_maps}(f), cannot be considered approximately single-valued and is therefore not useful as a representation of the R\"ossler flow.  
The return map construction fails because the 
variable $r$ does not vary monotonically for the points of intersection 
of the flow with $\poinc_2$, \cf \reffig{f:rossler_maps}(f). 
This difficulty motivates the use of manifold learning for the parameterization of the Poincar\'e 
section.

\subsection{\label{s:ML}Manifold learning}

Manifold learning or nonlinear dimensionality reduction~\cite{LeeVer07,Izenman12} is a sub-field of unsupervised machine learning.  The main problem it aims to address is projecting correlated data from a $d$-dimensional ambient space to a lower dimensional embedding space with dimension $d_r\ll d$, while preserving as much of the structure in the data as possible. In other words, if the data are thought to lie on (or close to) a manifold embedded in the high-dimensional space, then manifold learning aims to re-embed this manifold to a space of as low a dimension as possible, while preserving either the geometrical or topological properties of that manifold. In contrast to linear dimensionality methods, such as principal component analysis, manifold learning methods utilize nonlinear transformations, which enables them to re-embed manifolds generated by nonlinear processes.  A wide range of methods of nonlinear dimensionality reduction have been devised. Well known schemes which are of relevance to the study of nonlinear dynamical systems are Isomap\rf{tenenbaum2000}, which aims to devise transformations which preserve geodesic distance on the manifold, diffusion maps\rf{Coifman05,CoiLaf06,Nadler06}, which are robust even in the presence of noise, and locally linear embedding (LLE)\rf{RoSa00}, which is a topology preserving method. Here we choose to use LLE because it is relatively efficient, a well documented implementation exists in the \texttt{python} package \texttt{scikit-learn}\rf{scikit-learn} and topology (rather than distance) preservation is sufficient for our purposes. Similar results were also obtained for most cases studied here using Isomap.

LLE utilizes conformal mapping, which preserves local angles, in order to map data lying on (or close to) a manifold in a $d$-dimensional ambient space to the $d_r$-dimensional embedding space. The first step in the method is to select, for each data point, a neighborhood in which the manifold can be locally approximated as being linear. In practice this is achieved by detecting the $K$-closest neighbors of each data point, where $K$ an integer threshold. This allows each data point to be reconstructed from its neighbors, up to a certain reconstruction error. The key observation on which LLE is based is that the coefficients (or ``weights'') of the linear transformation  used for the reconstruction are (under suitable constraints) invariant under rotations, translations and rescalings. Thus, these coefficients are expected to be identical for both the original and the low-dimensional embedding of the data. By minimizing the global reconstruction error, one obtains global intrinsic coordinates which parametrize the manifold in the low dimensional space.  In practice, the problem of global error minimization is written in the form of a sparse eigenproblem of dimension $N_p\times N_p$, where $N_p$ is the number of data points. The coordinates of the data points in the $\dr$-dimensional re-embedding space are then identified with the eigenvectors corresponding to the smallest $\dr$ eigenvalues. 

There are two free parameters in LLE: the number of neighbors kept $K$ and a small normalization parameter $\delta$. The latter is required in order to condition the eigenvalue problem when the number of neighbors $K$ is larger than the number of ambient dimensions. 
Here, the parameter $K$ is chosen by estimating the number of neighbors within an approximately linear neighborhood in the Poincar\'e section.
It is noted that LLE can be formulated so that points within a ball of radius $\epsilon$, rather than $K$-nearest neighbors, are kept as part of the linear neighborhood. While that formulation would be more well suited for the study of chaotic dynamics, for which data are not uniformly sampled, it has not been used here since it is not implemented in the present version of \texttt{scikit-learn}. 
The normalization parameter $\delta$ was also varied and the construction of return maps repeated, in order to ensure that it leads to topologically equivalent results (by comparing kneading sequences, see \refsect{s:manif-poinc}).
We note that there exist variants of LLE\rf{DoGr03,ZhWa07}, which bypass the need for a normalization parameter at a slight increase in computational cost.

The computational cost of LLE (as implemented in \texttt{scikit-learn}) scales as $O[d \log(K) N_p \log(N_p)]+O[d N_p K^3]+O[\dr N_p^2]$,
where the first term refers to the detection of $K$ nearest neighbors, the second to the computation of local reconstruction weights and the third to the solution of the eigenvalue problem with dense methods\rf{scikit-learn-manifold}. In \texttt{scikit-learn} the cost associated to the last term is significantly reduced through the use of sparse eigensolvers. The scaling with respect to the number of data points $N_p$  makes it clear that it is much more efficient to first reduce the number of data points by introducing a Poincar\'e section and then perform LLE, rather than performing these operation in the reverse order.

\subsection{\label{s:manif-poinc}Manifold learning of Poincar\'e sections}

\begin{figure}[ht!]
    \includegraphics[width=\columnwidth]{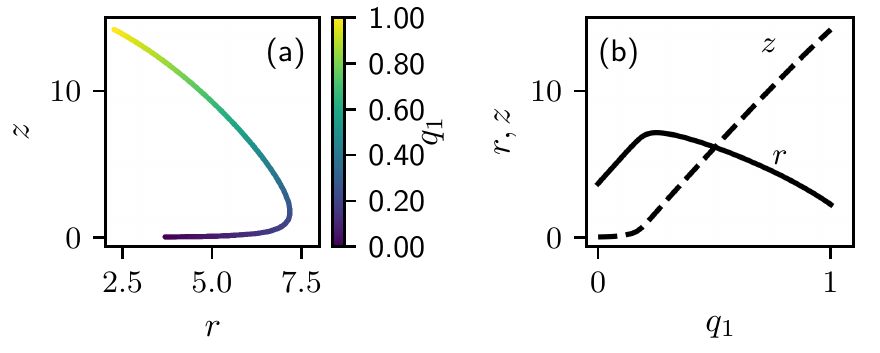}
        \caption{\label{f:rossler_s} (a) The points of intersection of the R\"ossler flow with the Poincar\'e section $\poinc_2$
        parametrized by the intrinsic variable $q_1$ obtained by LLE. (b) Interpolating functions $r(q_1)$ and $z(q_1)$ allowing to inverse LLE in order to
        obtain points in the original state-space for any value of $q_1$.
        }
\end{figure}

In the following we consider dynamical systems of the form
\beq\label{eq:ds}
    \dot{x}=F(x)\,,
\eeq
where $x\in\mathbb{R}^d$.
We use as input (training) data $\mathcal{D}$ the $N_p$ points of intersection of the attractor of \refeq{eq:ds} with the $(d-1)$-dimensional Poincar\'e surface of section. 
$\mathcal{D}$ is a fractal set of non-integer dimension $\df-1$, where $\df$ is the fractal dimension of the attractor. Our aim is to use manifold learning in order to re-embed the data to a space of integer dimension $\dr$, providing intrinsic variables $q\in \mathbb{R}^{\dr}$. The embedding dimension $d_r$ should be sufficiently high to allow the construction of a return map that captures the topology of the chaotic attractor. 
If some estimate of the attractor fractal dimension $\df$ is known, then the target dimension of the embedding space would be $\dr=\lfloor\df-1\rfloor$, where $\lfloor\,\rfloor$ is the floor function. 
However, one does not need to rely on estimates of the attractor fractal dimension but could instead increase $\dr$ iteratively until a single-valued return map is obtained, see \refsect{s:ks_2D}. Here we will restrict attention to cases for which we $\dr=1$ or $2$ is sufficient, but the same strategy is applicable for $\dr>2$.

As a first application we consider the points of intersection of the R\"ossler 
flow with the Poincar\'e sections $\poinc_1$ and $\poinc_2$. 
Application of LLE for both cases allows us to approximate the 
data of \reffig{f:rossler_maps}(b,e) as one-dimensional manifolds parametrized by a coordinate $q$, see \reffig{f:rossler_s}. 
We can thus construct 
the one-dimensional return maps of \reffig{f:rossler_maps}(d,f), which are single valued even for the ``problematic''
section $\poinc_2$ (up to the thickness discussed in \refsect{s:rossler}). The return maps of \reffig{f:rossler_maps}(c,d,f) describe the same dynamics and are thus conjugate to each other. We can show that this is indeed the case by examining the kneading sequence $K_c$, 
defined as the itinerary of the critical point, since two maps with the same kneading sequence are conjugate to each other~\cite{gilmore2003}.
For all three maps of \reffig{f:rossler_maps}(c,d,f) we find the kneading sequences $K_c=10010111111101\ldots$ to agree up to length at least 30.  

We will often need to iterate the return maps constructed through the use of LLE for out-of-sample data. This requirement will be handled by either linear or nearest-neighbor interpolation between sample points. An example where this is useful is when locating period-$p$ periodic orbits of the map by solving the condition $q_{1,n+p}=q_{1,n}$, through Newton–Raphson's method. 

There will also be occasions in which we need to extend the mapping $\mathcal{M}: x\mapsto q$, which relates in-sample data points to the LLE coordinates, to out out-of-sample points. This would allow to use the return map to predict the evolution of state-space points that are not part of the training data. In \texttt{scikit-learn} this can be achieved without re-computing the full transform through a mapping procedure that computes reconstruction weights only for the out-of-sample points, as described in\rf{SaRo02}.  A similar procedure can also be applied in order to construct an inverse LLE map $\mathcal{I}:\, q\mapsto x$, which is  however not implemented in \texttt{scikit-learn}. Here we take a simpler approach based on the fact that all return maps in this work can be eventually parametrized by a single LLE coordinate $q_1$ (\cf~\refsect{s:KSe}). Thus, it is sufficient to construct an inverse map as a set of $d$ interpolating functions $\mathcal{I}_j:\, q_1\mapsto x_j$, $j=1,\ldots,d$ using either nearest-neighbor or linear interpolation. \refFig{f:rossler_s}(a) and (b) illustrates the forward $\mathcal{M}: (r,z)\mapsto q_1$ and inverse $\mathcal{I}_1:\,q_1\mapsto r$, $\mathcal{I}_2:\,q_1\mapsto z$ maps, respectively, for the case of the R\"ossler example.

To illustrate the utility of the inverse map $\mathcal{I}$, suppose that we have obtained a root $q_1^*$  of the condition $q_{1,n+1}=q_{1,n}$ for the return map of \reffig{f:rossler_maps}(d). This can then mapped to a point $r^*,z^*$ on the Poincar\'e section by $\mathcal{I}_j$, $j=1,2$. This point is used as initial condition for \refeq{eq:rossler} and  integrated until the next intersection with the Poincar\'e section. The resulting solution is an approximation of a periodic orbit within plotting accuracy (that could be refined to desired accuracy by an iterative method~\cite{DasBuch}), as shown in \reffig{f:rossler_s}(a). A more systematic procedure for the determination of periodic orbits up to a given period is provided in \refapp{s:ks_cycles_method} and used in \refsect{s:KSe}.

Yet, one may argue that the R\"ossler example is somewhat artificial: \reffig{f:rossler_maps}(e) suggests that a usable return map could easily be constructed even for the case of section $\poinc_{2}$ by simply choosing $z$ as the parameterizing variable. 
For this reason we now turn to the discussion of application of manifold learning for Poincar\'e sections of a higher-dimensional system, 
for which simple parametrizations of the Poincar\'e section are not obvious, while the gains in dimensionality reduction will also be more significant.           

\section{\label{s:KSe}Application to Kuramoto-Sivashinsky equation}

\subsection{Kuramoto-Sivashinsky equation}

The Kuramoto-Sivashinsky equation (KSe)\rf{KurTsu76,siv}, 
\beq\label{eq:KSe}
    \partial_t u = (u^2)_x - u_{xx}-\nu u_{xxxx}\,,
\eeq
has received much attention as one of the simplest partial-differential equations (PDE) that exhibit spatiotemporally chaotic behavior. In \refeq{eq:KSe}, $u(x,t)$ is a one-dimensional scalar field, defined over a domain $x\in[-L/2,L/2]$ and $\nu$ is a damping parameter. Here we set $L=2\,\pi$ and use $\nu$ as control parameter. An alternative choice common in the literature is to set $\nu=1$ and use $L$ as control parameter. One can convert between the two conventions through~\rf{lanCvit07} $\nu=(2\,\pi/L)^2$. Here we follow \refref{Christiansen97} in imposing the Dirichlet boundary condition $u(-L/2,t)=u(L/2,t)=0$, which restricts
solutions to the subspace of antisymmetric functions, $u(x,t)=-u(-x,t)$.
Fourier expanding
\beq\label{eq:Fourier}
    u(x,t)=\sum_{k=-\infty}^{+\infty} \mathrm{e}^{\mathrm{i} kx}\,,
\eeq
taking into account the antisymmetry of $u(x,t)$, which implies $b_k=\mathrm{i}a_k$ with $a_k$ real, and  substituting in \refeq{eq:KSe} we obtain the infinite dimensional dynamical system 
\beq\label{eq:KSeFinf}
    \dot{a}_k= (k^2-\nu\, k^4)a_k - k \sum_{k=-\infty}^{+\infty} a_m\,a_{k-m}\,.
\eeq

In numerical simulations \refeq{eq:KSeFinf} is truncated to a finite dimensional set of ordinary differential equations by setting $a_k=0$ for $k>d$. Since $u(x,t)$ is real we have $a_k=-a_{-k}$ and choosing to work with zero-mean solutions, $a_0=0$, we obtain
 \begin{align}\label{eq:Fcoef_Trunc}
  \dot{a}_k & =  \left(k^2- \nu k^4\right)a_k  \nonumber \\
            &   - k\left( \sum_{m=1}^{k-1}a_m a_{k-m}-\sum_{m=k+1}^{d}a_m a_{m-k}
                    -\sum_{m=1}^{d-k}a_m a_{k+m}\right) \, ,
 \end{align}
where $k=1,\ldots,d$. Note that with our conventions the real field $u(x,t)$ is sampled on $2(d+1)$ points, since $a_{d+1}=a_0=0$. For the range of parameters studied here $d=16$ is sufficient for well resolved simulations~\cite{Christiansen97}. In the following we will study attractors of this 16-dimensional truncation of KSe. 

In the antisymmetric subspace, the continuous translational symmetry of KSe reduces to the discrete group of half-cell translations, $\tau u(x,t)=u(x+\pi,t)$\rf{KNSks90,SCD07}. This implies that if $u(x,t)$ is a solution then $u(x+\pi,t)$ is also a solution. Therefore solutions of \refeq{eq:KSe} either come in $\tau$-related pairs or are invariant (either pointwise or set-wise) under $\tau$~\rf{golubitsky2002sp}.  In terms of the Fourier modes half-cell translation symmetry takes the form 
\begin{align}\label{eq:translF}
    \tau a_{2m} &=a_{2m}\,,\\
    \tau a_{2m+1} &=-a_{2m+1}\,,
\end{align}
 where $m$ is an integer. The study of the system can be greatly simplified if points in state space related by symmetry are identified, \ie, if we apply symmetry reduction~\cite{DasBuch,GL-Gil07b}. This can be achieved by either a suitable change of variables (e.g. to invariant polynomials~\cite{GL-Gil07b}) or by introducing a fundamental domain~\cite{CvitaEckardt}, a subspace of the full state space which contains exactly one copy of each symmetry related pair. The latter approach is sufficient for our purposes and can be most easily achieved after the introduction of a Poincar\'e section. However, care must be taken that the Poincar\'e section is a group-invariant set, so that no state-space points leave the section under a symmetry transformation. This can be achieved by imposing a group-invariant condition to define the Poincar\'e section~\cite{siminosThesis}. Here we use the Poincar\'e section
 \beq\label{eq:KSePoinc}
    a_2 = 0\,,
 \eeq
subject to the orientation condition $\dot{a}_2>0$. \refEq{eq:KSePoinc} is invariant under the action \refeq{eq:translF} of $\tau$. An example of using this surface of section in a KSe simualtion with $\nu=0.0299$ is shown in \reffig{f:ks_1D}(a). We observe that the attractor consists of two distinct parts related by the action \refeq{eq:translF} of $\tau$. This motivates the definition of the fundamental domain consisting of all points $a$ with
\beq\label{eq:FD}
    a_3\geq0\,.
\eeq
 State space points that do not satisfy \refeq{eq:FD} are mapped to the fundamental domain by applying the half-cell translation operation $\tilde{a}=\tau a$, where $\tau$ is given by \refeq{eq:translF} and $\tilde{a}$ denotes a point in the fundamental domain. This leads to the symmetry reduced Poincar\'e section of \reffig{f:ks_1D}(b). In what follows we will simplify notation by referring to both points in the full space and the fundamental domain by the symbol $a$.  
We note that if we have chosen periodic boundary conditions then the symmetry group of KSe would be $O(2)$ and we would need to also implement continuous symmetry reduction~\cite{rowley_reconstruction_2000,BeTh04,SiCvi10,BudCvi14}. 
 
In practice it is more efficient to use a pseudospectral evaluation\rf{boyd01} of the nonlinear term in \refeq{eq:Fourier}, rather than the explicit sum in \refeq{eq:Fcoef_Trunc}. Due to stiffness in KSe we use an implicit scheme based on a backward differentiation formula~\cite{ShaRei97}, as implement in the library \texttt{scipy}\rf{2020SciPy-NMeth}. Variable time-step is used with absolute and relative tolerance set to $10^{-9}$ and $10^{-6}$, respectively. It was checked that increasing these tolerances by an order of magnitude results in the same kneading sequences up to topological length at least 10 for the return maps presented here.

Return maps for KSe with Dirichlet boundary conditions have been studied in \refref{Christiansen97,lanCvit07}. In \refref{Christiansen97} a one-dimensional return map is constructed for $\nu=0.029910$, using the geodesic distance measured along an ordered set of points of periodic points on the Poincar\'e section. However, in this case a single Fourier mode would also be sufficient for the parameterization of the return map. For a smaller value of the damping parameter, $\nu=0.026634$, \refref{lanCvit07} showed that there is considerably more complex structure in the attractor.  The latter was decomposed in two parts, each parametrized by the geodesic distance along the one-dimensional unstable manifolds of two short periodic orbits, resulting in two disjoint one-dimensional return maps. 

Here we choose to work with values of the damping parameter between these two extremes. The first case, with $\nu=0.0299$ has a Poincar\'e section for which no obvious parameterization by a single Fourier mode could be found, yet can still be well approximated as one-dimensional using manifold learning. For the second case we choose $\nu=0.02973$ which corresponds to one of the simplest cases for which a two-dimensional embedding is required.

\subsection{\label{s:ks_1D}One-dimensional embedding}

\begin{figure*}[ht!]
     \includegraphics[width=\textwidth]{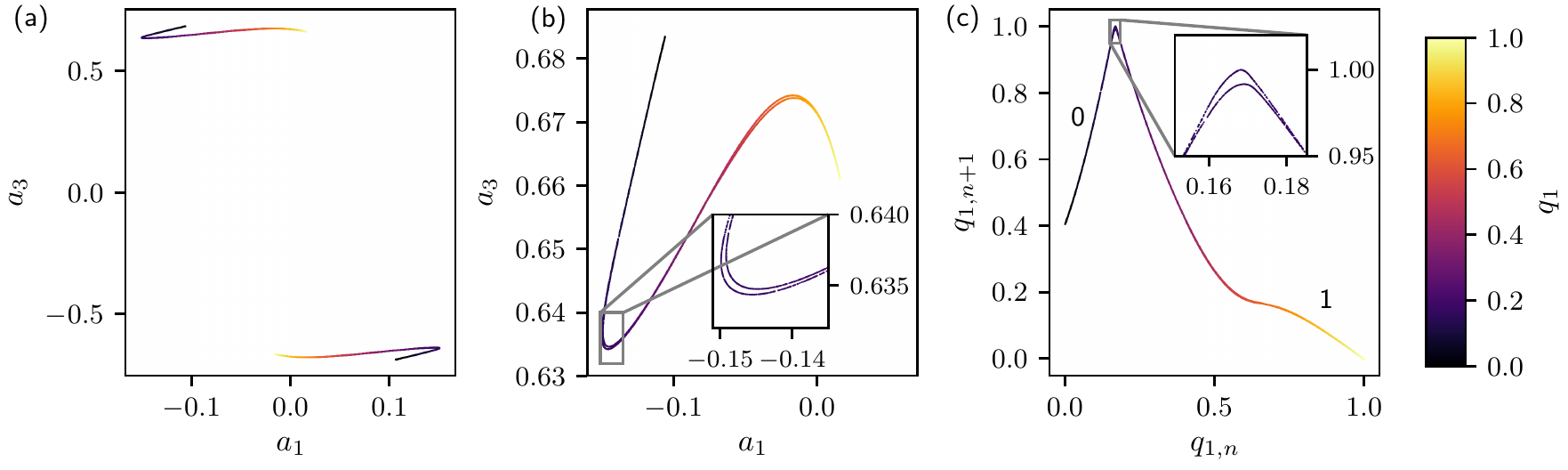}
         \caption{\label{f:ks_1D}
         KSe dynamics for $\nu=0.0299$. (a) Full space Poincar\'e section, $a_1-a_3$ modes projection, (b) Poincar\'e section reduced to the fundamental domain defined by $a_3\geq0$. 
         (c) Return map using a 1D LLE embeddeding, with $K=48$ nearest neighbors and regularization parameter $\epsilon=10^{-3}$.
        }
\end{figure*}

For the first of the cases that we study we take $\nu=0.0299$ and integrate KSe up to $t=4000$, which results in $N_p=9092$ intersection points with the Poincar\'e section. A projection of the Poincare section in the full space and fundamental domain is shown in \reffig{f:ks_1D}(a) and \reffig{f:ks_1D}(b), respectively. The footprint of the attractor on the Poincar\'e section is extremely thin, suggesting that a one-dimensional embedding should be sufficient to represent the dynamics. However, \reffig{f:ks_1D}(b) reveals no obvious variable for the construction of a return map. Manifold learning with LLE solves this problem in a natural way, allowing to construct a 1D embedding, visualized by the color-code used in \reffig{f:ks_1D}(b). This results to a unimodal return map $q_n\mapsto q_{n+1}$ on the unit interval, shown in \reffig{f:ks_1D}(c).
The critical point $x_c$ partitions the unit interval to two base intervals which we label with the symbols $0$ and $1$. The kneading sequence associated with the critical point is $K_c=1011111111\ldots$ and will be used in \refsect{s:ks_cycles} to determine the spectrum of admissible periodic orbits.

We note that this map is only approximately unimodal, since the intersection of the attractor with the Poincar\'e section is only approximately one-dimensional, see the insets in \reffig{f:ks_1D}(c) and (b), respectively. We see in particular that there is transverse structure in the return map of \reffig{f:ks_1D}(c), where we can discern an ``upper'' or ``primary'' branch and a ``lower'' or ``secondary'' branch. The origin of these branches will be discussed in \refsect{s:ks_UM}. Zooming in further would reveal additional structure (not shown). We could try to refine our description by keeping a second dimension in LLE, as we do for the lower damping parameter case in \refsect{s:ks_2D}. However, this will not be necessary for the present example, since we will show in \refsect{s:ks_cycles} that the approximate one-dimensional map correctly encodes topological information on the spectrum of short periodic orbits of the system.

\subsection{\label{s:ks_2D}Two-dimensional embedding case: from a tree map to a trimodal map}

\begin{figure*}[ht!]
    \includegraphics[width=\textwidth]{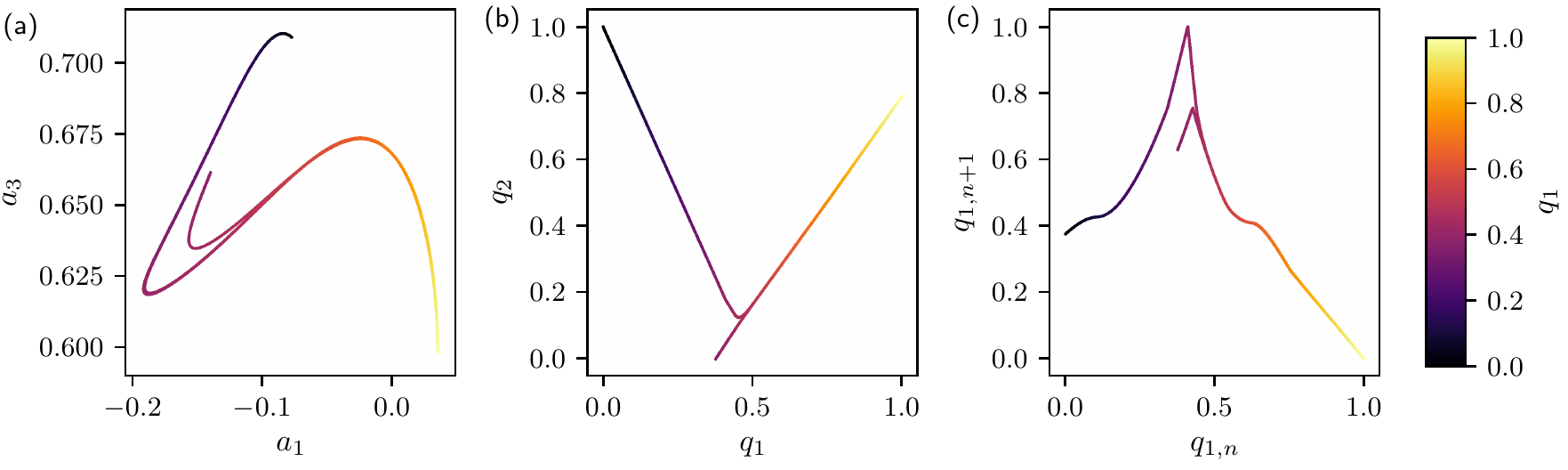}
        \caption{\label{f:ks_2D}
        (a) Projection of KSe Poincar\'e section on $a_1,\,a_3$ modes for $\nu=0.02973$ (in the fundamental domain). (b) Two dimensional embedding of the Poincar\'e section on LLE coordinates $(q_1, q_2)$, using $K=143$ nearest neighbors and regularization parameter $\epsilon=0.5\times10^{-3}$.
        (c) Naive return map using the coordinate $q_1$ is multi-valued.
        }
\end{figure*}

For slightly smaller value of the damping parameter, $\nu=0.02973$, KSe was integrated up to $t=40000$, resulting in $N_p=93640$ intersection points with the Poincar\'e section, shown in \reffig{f:ks_2D}(a).
A naive return map $q_{1,n}\mapsto q_{1,n+1}$, \reffig{f:ks_2D}(c), using the first LLE embedding variable $q_1$ results in a return map which cannot be considered approximately single valued, \ie, we cannot ignore the secondary branch of the map. It should be noted that such return maps can still be studied if one augments them with grammar rules specifying whether the ``upper'' or ``lower'' branch of the map would be visited in the next intersection~\cite{hansen}. However, this complicates the analysis of symbolic dynamics and we will instead develop a methodology that to would allow to disentangle maps like the one in \reffig{f:ks_2D}(c). The origin of the multivaluedness of this map can be traced back to \reffig{f:ks_2D}(a), where we can see that due to the existence of a secondary branch in the Poincar\'e section, $q_1$ cannot uniquely parameterize all points on the section. A two dimensional embedding is thus computed  through application of LLE on the Poincar\'e section data, as shown in \reffig{f:ks_2D}(b), where $q_1$ and $q_2$ are the embedding variables. This would seem to imply that we need to construct a two-dimensional return map $(q_{1,n}, q_{2,n}) \mapsto (q_{1,n+1}, q_{1,n+1})$. However, an alternative route is taken here, which will allow us to construct a one-dimensional map on the unit interval with four-letter alphabet.

\begin{figure}[ht!]
    \includegraphics[width=\columnwidth]{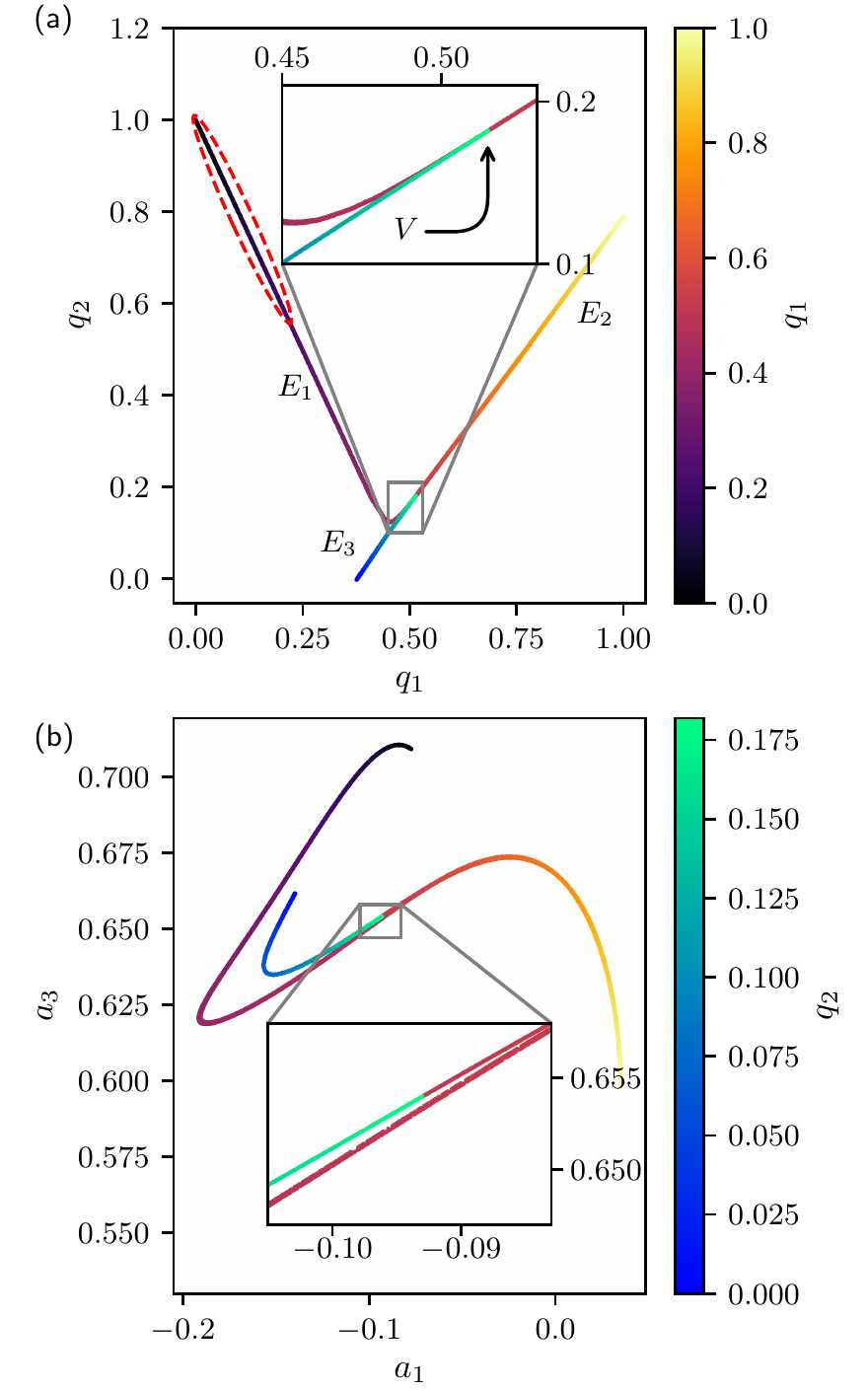}
        \caption{\label{f:ks_split}(a) Two-dimensional LLE embedding of the Poincar\'e section of \reffig{f:ks_2D}(b) interpreted as a tree with edges $E_i$, $i=1,\,\ldots,\,3$, meeting at the vertex $V$ (see inset). Points on egdes $E_1$ and $E_2$ are parametrized (and color-coded) by the first LLE coordinate $q_1$, while points on $E_3$ are parametrized (and color-coded) by the second coordinate $q_2$. The red, dashed ellipse encloses the subset $E_s$ of $E_1$ consisting of points mapped to $E_3$ in the next intersection with the Poincar\'e section. (b) Projection of KSe Poincar\'e section on $a_1,\,a_3$ modes, following the same color-coding as in panel (a).
        }
\end{figure}

Our starting point is the observation that the embedding of \reffig{f:ks_2D}(b) can be thought of as (approximately) having the topology of a tree 
with three edges, $E_i$, $i=1,\ldots,3$, and an internal vertex $V$, see \reffig{f:ks_split}(a). Intuitively, we observe that edges $E_1$ and $E_2$ can be parametrized by $q_1$ while edge $E_3$ can be parametrized by $q_2$.
However, in order to realize such a parameterization of the computational data we need a way to uniquely assign data points to each edge. The availability of a two-dimensional LLE embedding space is particularly invaluable in this respect, since it makes the problem of assigning data to the different edges much more tractable than in the original 15-dimensional Poincar\'e section.
To this end, we perform a linear fit of the data shown in the inset of \reffig{f:ks_split}(a) to obtain a curve $\mathcal{C}$ which separates points on $E_3$ from points on $E_1$ and $E_2$. 
Then a data point belongs to $E_3$ if it satisfies 
$q_2<\xi(q_1)$,
subject to the additional condition $q_1<q_{1,c}$, or belongs to $E_1\cup E_2$ otherwise. Here, 
$\xi(q_1)$ is the graph of $\mathcal{C}$ and $q_{1,c}$ effectively defines the position of the vertex $V$, as discussed shortly. 

\begin{figure*}[ht!]
    \includegraphics[width=\textwidth]{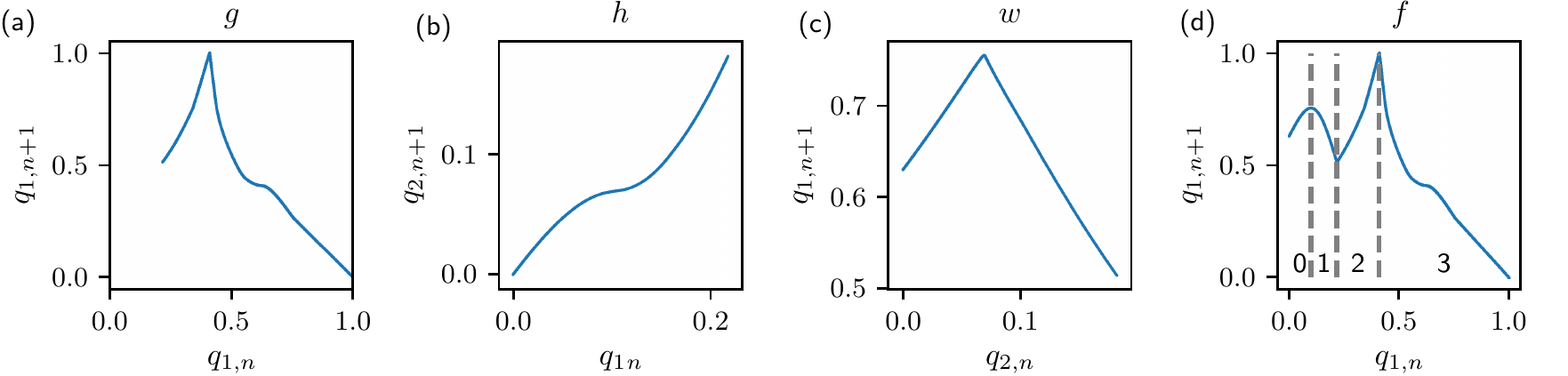}
        \caption{\label{f:ks_maps} Tree maps and combined map on the unit interval (a) map $g:\, E_p \mapsto E_1\cup E_2$, (b) $h:\, E_r\, \mapsto E_3$, (c) $w:\, E_3\, \mapsto E_2$,  (d) $f:\, E_1\cup E_2 \mapsto E_1\cup E_2$. The dashed lines indicate the separation of the invariant interval $[0,1]$ to sub-intervals corresponding to different symbols $s=0,\,\ldots,\,3$.
        }
\end{figure*}

With this construction we are now in position to identify which points on $E_1\cup E_2$ are mapped on $E_1\cup E_2$ and which on $E_3$ in the next intersection with the Poincar\'e section, and vice-versa. This allows to define three maps between edges of the tree:
\begin{enumerate}
 \item $g:\, E_p \mapsto E_1\cup E_2$, where $E_p\subset E_1\cup E_2$. In terms of \reffig{f:ks_split}(a) the set $E_p$ corresponds to all points in $E_1\cup E_2$ except those encircled by the red ellipse and is pararametrized by $q_1$. The map $g: q_{1,n}\mapsto q_{1,n+1}$ is shown in \reffig{f:ks_maps}(a). 
 \item $h:\, E_s\, \mapsto E_3$, where $E_s\subset E_1$ is encircled by the red ellipse in \reffig{f:ks_2D}. 
 The map $g: q_{1,n}\mapsto q_{2,n+1}$ is shown in \reffig{f:ks_maps}(b).
 \item $w:\, E_3\, \mapsto E_r\subset E_2$. The graph of $w:q_{2,n}\mapsto q_{1,n+1}$ is shown \reffig{f:ks_maps}(c). As will be discussed shortly, the restriction of the range of w to a subset of $E_2$ (rather than a subset of $E_1\cup E_2$) can been achieved by a suitable choice of the position of the vertex $V$, through the parameter $q_{1,c}$, which defines the boundary between $E_3$ and $E_2$. This choice allows to ensure continuity of the the combined map $f$, which will be constructed in the following. 
\end{enumerate}

We have thus constructed a description of the KSe attractor in terms of three coupled repellers approximated as maps on a tree. Such tree maps have appeared in the study of pruning of the H\'{e}non map\rf{CaHa02}. Since the kneading theory for such maps is quite involved~\cite{AlvRa04} a more direct approach is introduced here, which consists in using function composition in order to construct a single return map on the unit interval. In particular, we define the map:
\beq
f(q_1)= 
    \begin{cases}
        w\left(h(q_1)\right)\,,& \text{if } q_1< q_{1,c}\\
        g(q_1),              & \text{if } q_1\geq q_{1,c}\,,
    \end{cases}
\eeq
shown in \reffig{f:ks_maps}(d). We have thus reduced the dynamics on the two-dimensional Poincar\'e section to a one-dimensional trimodal (\ie, having three maxima and minima) map on the interval. This map is approximately single-valued and continuous (up to fine structure similar to that discussed in \refsect{s:ks_1D}). The continuity property can be imposed by a proper choice of the position of the vertex $V$. Since we deal here with a finite data set, we have some freedom to choose the position of this vertex, effectively changing the assignment of data points to different edges. 
We therefore choose a cutoff value $q_{1,c}$ which results in a continuous map, in the sense that any discontinuity in the map is of the order of the average distance between data points in our sample.
If we had chosen a
different a value $q_{1,c}'<q_{1,c}$ or $q_{1,c}'>q_{1,c}$, then we would have $w(\left(h(q_{1,c}')\right)>g(q_{1,c}')$ or
$w(\left(h(q_{1,c}')\right)<g(q_{1,c}')$, respectively, and the map $f$ would be discontinuous. While we could still obtain symbolic dynamics using a discontinuous map\rf{SeCo10}, our study is significantly simplified by ensuring continuity of $f$ within the limits set by the finite sample size.
Note that by this choice the tree topology is ensured for the two-dimensional LLE embedding but not necessarily for the the original 15-dimensional Poincar\'e section, compare the insets in \reffig{f:ks_split}(a) and \reffig{f:ks_split}(b). As we will show, this choice to model the dynamics based on the LLE embedding rather than the original space correctly captures the spectrum of short (up to topological length $n=10$) periodic orbits of the system, see \refsect{s:ks_cycles}.
Note that the map $f$ is not a first return map to the Poincar\'e section, but could be thought of as a first return map from the subset $E_1\cup E_2$ to itself. 
 
The three critical points $q_{i,c}$, $i=1,\ldots,3$ partition the unit interval into four sub-intervals which can be assinged the symbols 0 to 4, as shown in \reffig{f:ks_maps}(d). There are three corresponding kneading sequences, $K_1=3233303232\ldots$, $K_2=3333333323\ldots$ and $K_3=3032303230\ldots$, which encode the topological organization of the flow and determine the spectrum of admissible periodic orbits, see \refsect{s:ks_cycles}. It was verified that the kneading sequences do not change up to topological length 10 when either the numerical integration tolerance parameters are reduced by an order of magnitude, integration time was increased by 50\%, or LLE parameters were varied.

\subsection{\label{s:ks_cycles}Periodic orbits}

One of the advantages of representing continuous time dynamical systems with discrete time return maps, is that the latter allow the systematic 
determination of unstable periodic orbits which provide a skeleton of the dynamics and could be used in cycle expansions in order to compute average values of observables~\cite{DasBuch}. For KSe in the antisymmetric domain there are two types of 
such solutions: periodic orbits with period $T_p$, satisfying $u(x,T_p+t)=u(x,t)$, and pre-periodic (or relative periodic) orbits $u(x,T_p+t)=\tau u(x,t)$, which are periodic of period $T'_p=2\,T_p$ in the full space and of period $T_p$ in the fundamental domain. We will refer to both types of orbits as \emph{cycles} and distinguish them by their multiplicity in the full space: periodic orbits come in pairs related by half-cell translation $\tau$ (multiplicity $m=2$), while pre-periodic orbits are $\tau$-invariant sets (multiplicity $m=1$).

For each fundamental domain cycle we can define its topological length $n$ as the number of intersections of the orbit with the Poincar\'e section (restricted to $E_1\cup E_2$ for the tree-topology case, $\nu=0.02973$).  Note that a cycle with multiplicity $m=1$ and topological length $n$ in the fundamental domain will have topological length $2n$ in the full space. If our discrete time return maps are indeed a valid representation of the continuous time flow
then there must be a 1-1 correspondence between the return map cycles of period $n$ and the flow cycles of topological length $n$. This allows us to both detect 
the cycles of the flow in a systematic manner and to test the limits of validity of the 1-D map approximation.

The procedure we follow in order to systematically detect cycles up to a given topological length is based on determining the admissible cycles through the well known kneading theory of one-dimensional maps~\cite{MilThu88}, see\rf{DasBuch,gilmore2003,hansen} for physicist-friendly expositions. In order to formulate guess periodic orbit solutions for an iterative solver we exploit the partitioning of the invariant interval to sub-intervals by the pre-images of the kneading sequence. Since this procedure has not, to the author's knowledge, been reported elsewhere it is briefly outlined in \refapp{s:ks_cycles_method}.

\begin{table}[ht]
\caption{\label{t:cycles}All cycles up to topological length $6$ for KSe, damping parameter $\nu=0.0299$ (unimodal map) and $\nu=0.02973$ (trimodal map). 
Listed are the cycle label $p$ (overbars are dropped), fundamental domain period $T_p$, fundamental domain leading Floquet exponent $\Lambda_1$ and full-space cycle multiplicity $m$.}
    \begin{ruledtabular}
        \begin{tabular}{rrrr}
        {$p$} & {$T_p$} & {$\Lambda_1$} & {$m$}\\
        \hline
         \multicolumn{4}{l}{Unimodal, $\nu=0.0299$} \\
            1             & 0.449 & -1.838  &  1 \\
            01            & 0.870 & -2.085  &  2 \\
            0111          & 1.751 & -4.206  &  2 \\
            010111        & 2.630 & 7.800   &  2 \\ 
            011111        & 2.640 & -7.132  &  2 \\
            \hline
         \multicolumn{4}{l}{Trimodal, $\nu=0.02973$} \\ 
            3       &  0.447  &-2.227     & 1\\ 
            23      &  0.857  &-3.230     & 2\\  
            1323    &  2.144  &-10.792    & 1\\ 
            1333    &  2.176  & 21.297    & 1\\  
            2333    &  1.737  &-12.468    & 2\\ 
            032323  &  2.980  &-12.143    & 1\\ 
            132323  &  2.995  & 21.210    & 1\\ 
            133323  &  3.043  &-79.969    & 1\\ 
            133333  &  3.072  & 115.592   & 1\\ 
            232333  &  2.599  & 41.663    & 2\\  
            233333  &  2.627  &-54.925    & 2
        \end{tabular}
    \end{ruledtabular}
\end{table}

For the KSe example with $\nu=0.0299$ all cycles up to topological length $n=8$ determined by kneading theory were found to correspond to cycles of the continuous time flow. Table \ref{t:cycles} lists some of the shortest cycles as well as their leading Floquet multiplier (for pre-periodic orbits cycles this is computed in the fundamental domain) and their multiplicity. At topological length $n=10$, the boundary value solver did not converge for two of the cycles permitted by kneading theory, $\overline{0101111111}$ and $\overline{0111111111}$. The reason for this is the transverse structure of the return map shown in the inset of \reffig{f:ks_1D}(c), which implies additional pruning rules due to the secondary (lower) branch of the map. Indeed, the right-most point of each of these two cycles corresponds to topological coordinate $\theta$ exceeding that of the image of the critical point associated with the secondary branch, implying that both cycles are pruned (\cf \refapp{s:ks_cycles_method}). In order to formulate such additional pruning rules in a consistent manner we would need to introduce a second LLE coordinate and eventually construct a multi-modal map, as we did for the case $\nu=0.02973$. However, this was not attempted since we are interested here in the coarsest level of organization of the flow. 

\begin{figure}[ht!]
    \includegraphics[width=0.6\columnwidth]{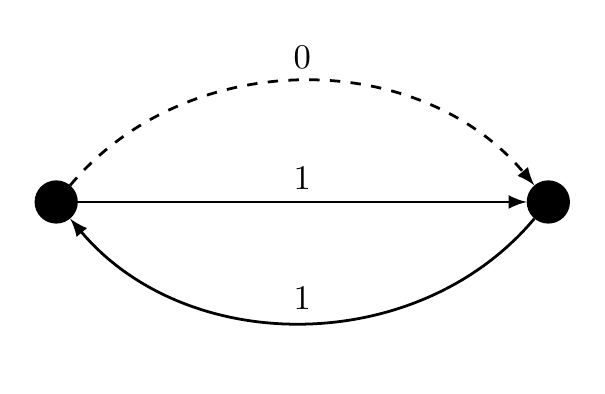}
        \caption{\label{f:automaton} Automaton representation of the finite grammar approximation of KSe recurrent set symbolic dynamics for $\nu=0.0299$. Traversing a dashed line corresponds to symbol $0$ while traversing a solid line corresponds to $1$.
        }
\end{figure}

A finite grammar approximation of the symbolic dynamics for $\nu=0.0299$ can be obtained 
by approximating the kneading sequence by $K_c=10\overline{1}$. This leads to two forbidden sequences (pruning blocks), $00$ and $0110$. The associated representation of the recurrent set symbolic dynamics as a finite automaton~\cite{hansen} is shown in \reffig{f:automaton}.  It correctly encodes all admissible periodic orbits up to length $n=8$, but would need to be supplemented with additional rules at length $n=10$. 
Using this graph the corresponding characteristic polynomial can be determined\rf{DasBuch} to be
\begin{equation}\label{eq:char_pol_ks_1D}
 p(z)=1-2z^2\,.
\end{equation}
It allows us to compute the topological entropy $h=-\ln(z_0)\simeq0.35$, where $z_0=1/\sqrt{2}$ solves $p(z_0)=0$. The growth of the number of admissible intervals with topological length $n$ found by admissibility criteria imposed by kneading theory (\cf \refapp{s:ks_cycles_method}) is in good agreement with the estimate $\sim e^{h\,n}$; however, this should be  considered an upper bound since additional pruning rules would need to be introduced for for $n>8$. 

\begin{figure*}[ht!]
    \includegraphics[width=\textwidth]{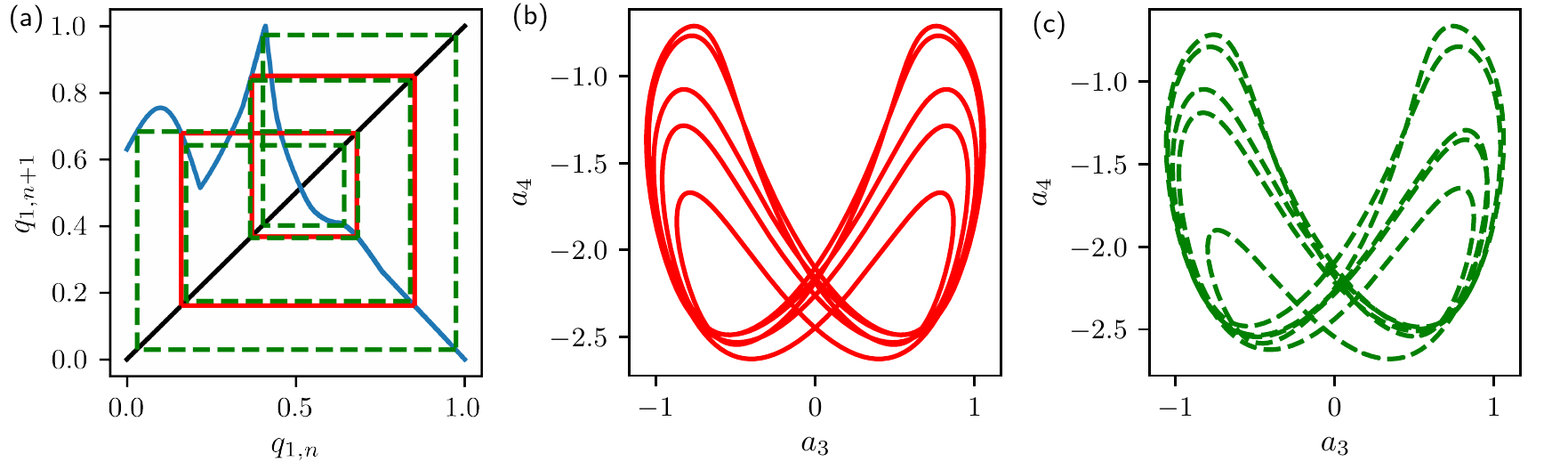}
        \caption{\label{f:ks_cycles} KSe cycles for $\nu=0.02973$. (a) Cycles $\overline{1323}$ (red, solid lines) and $\overline{03231323}$ (green, dashed lines) of the discrete time return map represented by 'cobwebs', (b,c)  cycles $\overline{1323}$ and $\overline{03231323}$, respectively, of the continuous time flow, projected on $a_3$, $a_4$ modes.
        }
\end{figure*}

For the multi-modal map case, $\nu=0.02973$, all admissible cycles of the discrete map up to length $n=10$ led to converged searches for continuous time flow cycles. The shortest cycles, up to length $n=6$ are listed in \reftab{t:cycles}, while two cycles are plotted in \reffig{f:ks_cycles}. In order to obtain reliable initial conditions for the detection of some of the cycles of length $n=10$ we had to increase the number of data points used in LLE (by increasing integration time) by 50\%. The periodic orbit data (as well as the kneading sequences) suggest that symbols $0$, $1$, $2$ and $3$ are always followed by $3$;
however there exist multiple additional grammar rules of various lengths, suggesting that a useful finite grammar approximation might be out of reach in this case.

\subsection{\label{s:ks_UM}Organization of the flow in terms of unstable manifolds}

\begin{figure}[hb!]
    \includegraphics[width=\columnwidth]{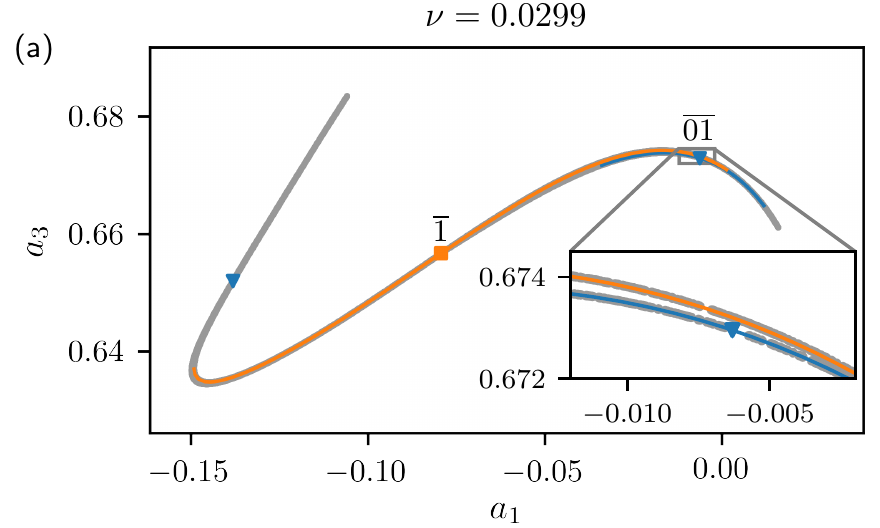}\vspace{10pt}\\
    \includegraphics[width=\columnwidth]{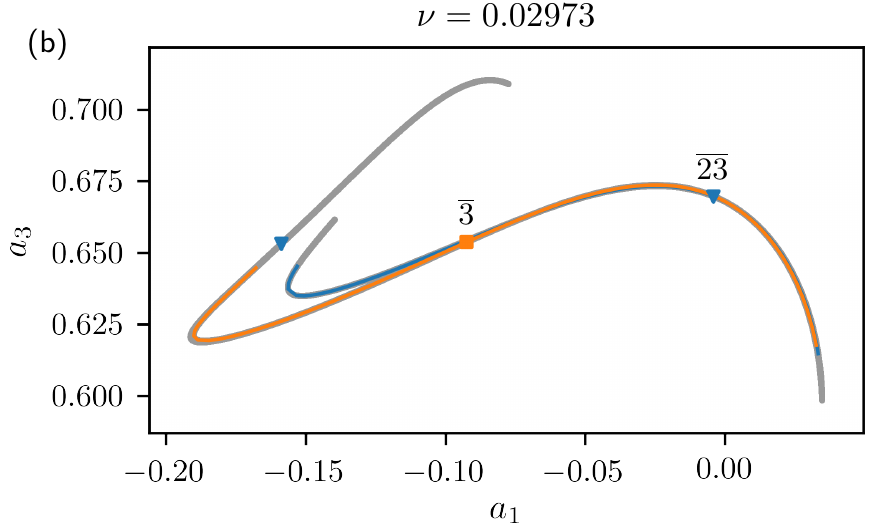}
        \caption{\label{f:ks_UM} (a) KSe Poincar\'e section projection for $\nu=0.0299$, showing a long trajectory (grey), periodic points of cycles $\overline{1}$ (orange square) and $\overline{01}$ (blue triangles) and corresponding parts of their unstable manifolds, (orange/blue solid lines, respectively). (b) The same for $\nu=0.02973$, showing periodic points of cycles $\overline{3}$ (orange square) and $\overline{23}$ (blue triangles) and corresponding parts of their unstable manifolds, (orange/blue solid lines, respectively). For cycles $\overline{01}$ and $\overline{23}$, only the part of the unstable manifold which connects to one of the periodic points has been visualized, for clarity.
        }
\end{figure}

One way to understand the topological organization of a chaotic flow is by tracking the stretching and folding of the unstable manifolds of compact solutions, such as periodic orbits~\cite{DasBuch}. Unstable manifolds stretch due to local instability until they are folded by nonlinearity, without self-intersections. The process repeats itself creating the ``transverse structure'' in the Poincar\'e section of \reffig{f:ks_1D}(b). 
An alternative, computationally more tractable way to understand this effect is by examining  finite segments of the unstable manifolds of a few of the shortest cycles embedded in the attractor, as we now show.

All KSe cycles detected in this work have a single unstable direction associated with a real Floquet eigenvalue. For the unimodal map case, $\nu=0.0299$,  we plot in \reffig{f:ks_UM}(a) segments of the global unstable manifolds $\mathcal{U}_{\overline{1}}$, $\mathcal{U}_{\overline{01}}$ of the two shortest cycles $\overline{1}$ and $\overline{01}$, respectively. These manifolds have been approximated by integrating forward in time $50$ trajectories with initial conditions on the local unstable manifold of each cycle. We choose the integration time to be long enough to capture a significant segment of the global manifolds, yet short enough to avoid the inevitable nonlinear ``folding'' of the manifolds.  As the inset in \reffig{f:ks_UM}(a) shows each of these unstable manifolds spans a different ``branch'' of the transverse structure of the attractor. Longer cycles could be used to approximate the transverse structure with finer detail.
The one-dimensional approximation amounts to disregarding the transverse structure. In this respect our treatment is similar to the parametrization of the the Poincar\'e section by the geodesic length along the unstable manifold of a short cycle proposed in \refref{lanCvit07}.

The trimodal map case, $\nu=0.02973$, presents particular interest.
 Segments of the global unstable manifolds $\mathcal{U}_{\overline{3}}$, $\mathcal{U}_{\overline{23}}$ of the two shortest cycles $\overline{3}$ and $\overline{23}$, respectively, are plotted in \reffig{f:ks_UM}(b). The unstable manifold $\mathcal{U}_{\overline{3}}$ spans much of the primary branch of the attractor (which we parameterized by $q_1$), while $\mathcal{U}_{\overline{23}}$ extents from the primary branch to the secondary branch (which we parameterized by $q_2$). 
The ``mis-alignment'' of these unstable manifolds leads to the need for a two-dimensional embedding. At the same time the fact that the unstable manifolds are one-dimensional explains why the tree approximation of the Poincar\'e section data was successful and led to the reduction of the dynamics to a one-dimensional map.

\section{\label{sec:conclusions} Conclusions}

In this work, it was shown through numerical examples that manifold learning can be used as a means to uncover the low-dimensional topological organization of high-dimensional chaotic dynamical systems by facilitating the construction of return maps. By contrast to previous work, nonlinear dimensionality reduction is applied directly on data on the Poincar\'e section, rather than on the continuous time attractor. This turns out to be computationally more efficient, as it signifcantly reduces the amount of data required, while preserving all information required for the construction of return maps. More importantly, it allows to decouple the problem of choosing a Poincar\'e section from the problem of parameterizing the return map, as discussed in \refsect{s:rossler}. By introducing an intrinsic parameterization of the data, manifold learning allows one to choose any Poincar\'e section satisfying the minimal requirement of being transverse to the flow. If no single Poincar\'e section satisfying this requirement could be found, local sections could be used~\cite{atlas12} and the resulting local maps combined as for the case of the tree maps of \refsect{s:ks_2D}. 

Application of the method on KSe allowed us to reduce the dynamics from the $d=16$-dimensional state-space in which well resolved solutions of the PDE live to one-dimensional maps. The latter enabled the systematic determination of the periodic orbits of the KSe flow, up to topological length $n=10$. For the KSe example with lower value of the damping parameter, $\nu=0.02973$, a two-dimensional embedding space for the Poincar\'e surface of section was required. Due to strong dissipation, the problem turned out to be quasi one-dimensional, with the data having the approximate topology of a tree. As a first step, a tree map was constructed which eventually was reduced  to a one-dimensional trimodal map. The reason that this was possible was traced back to the fact that all unstable manifolds of periodic orbits, which organized the attractor, were one-dimensional. We may thus conjecture that for a wide class of systems possesing compact solutions with a low number of unstable directions a similar reduction procedure should be possible. The factor determining the dimensionality of the resulting return maps would be the maximal number of unstable directions, rather than the embedding dimension of the attractor.

\begin{acknowledgments}
I would like to thank Alexander Jonsson, Caroline B\"ukk and Oscar Johansson for early contributions to this project
and Predrag Cvitanovi\'c for inspiring discussions.
This research has been supported by the Swedish Research Council, Grant No. 2016-05012.
\end{acknowledgments}

\section*{Data Availability Statement}

The data that support the findings of this study are available from the author upon reasonable request. The code used in this work is available as a library in the online repository \url{https://github.com/vasimos/poinc-man}.

\appendix

\section{\label{s:ks_cycles_method}Methodology for periodic orbit searches}

Here we outline the procedure developed for the systematic determination of periodic orbits of a continuous time flow, given a parameterization by an approximate one-dimensional map.

The first step is to determine the admissible cycles of topological length $n$. Consider a $M$-modal map $f$ of the unit interval, with critical points $x_{c,i}$, $i=1,\ldots,M$, separating the invariant interval, here $I=[0,1]$, to $N=M+1$ base intervals $I_k$, $k=0,\ldots,M$ to each of which we associate a symbol $s=0,\ldots,M$. To simplify notation assume that the branches of $f$ corresponding to even symbols are orientation preserving [$f'(x)>0$], while those corresponding to odd symbols are orientation reversing [$f'(x)<0$]. The algorithm can be easily adapted if this condition is not met. Take the following steps.
\begin{enumerate} 
 \item Determine the kneading sequences $K_i$ of all critical points $x_{c,i}$ of the map, as well as the associated kneading values $\theta_i=\theta(K_i)$,  where $\theta$ maps symbolic sequences to the so-called topological coordinate~\rf{gilmore2003}. The latter is defined for any symbolic sequence $\Sigma=s_0s_1\ldots$ as
 \[
  \theta(\Sigma) = \sum_{i=0}^{+\infty}\frac{t^i}{N^{i+1}}\,,
 \]
 where
 \[
   t_i=    
    \begin{cases}
        s_i\,,		& \text{if } \epsilon(s_0)\epsilon(s_1)\ldots\epsilon(s_{i-1})=1\,,\\
        (N-1)-s_i\,,	& \text{if } \epsilon(s_0)\epsilon(s_1)\ldots\epsilon(s_{i-1})=-1\,,
    \end{cases}
 \]
 and the parity operator $\epsilon(s_k)=1$ ($\epsilon(s_k)=-1$) if $s_k$ is even (odd).
 The topological coordinate $\theta\in[0,1]$ spatially orders symbolic sequences (both admissible and inadmissible ones).
 \item Generate all prime cycle labels of length $n$ using Duval's algorithm\rf{Duval88} to generate necklaces (periodic symbol sequences).
 \item For each prime cycle label generate the symbolic sequences of all $n$ cycle points by cyclic permutations of the cycle label.
 \item We can use the fact that critical points are local maxima or minima of the map $f$ in combination with spatial ordering of sequences by the topological coordinate $\theta$ to arrive at an admissibility condition for periodic orbits, see~\cite{hansen} for a readable review for multi-modal maps. For each critical point $x_{c,i}$ identify the set of all -points of a given cycle within the intervals $I_{i-1}$, $I_{i+1}$, i.e. to the left and right of the critical point. Denote by $s_0s_1\ldots s_{n-1}$ the symbolic sequence of any cycle point within this set and by $s_1\ldots s_{n-1}s_0$ its forward image. Denote by $\theta^{\mathrm{max}}(s_1\ldots s_{n-1}s_0)$ [$\theta^{\mathrm{min}}(s_1\ldots s_{n-1}s_0)$] the maximum [minimum] value of the topological coordinate over the set of the forward images. Then a cycle is inadmissible if and only if there exists a critical point $x_{c,i}$ for which
  \begin{align}
    \theta^{\mathrm{max}}(s_1\ldots s_{n-1}s_0)>K_i, & \quad\text{for i even}\,,\\
    \theta^{\mathrm{min}}(s_1\ldots s_{n-1}s_0)<K_i, & \quad\text{for i odd}\,.   
  \end{align}
\end{enumerate}

Having determined the set of admissible orbits we then need to generate initial guesses for the periodic points of the return map $f$. Our strategy is to determine, for each periodic point $s_0s_1\ldots s_n$ of the map, sub-intervals $I_{s_0s_1\ldots s_n}$ of the invariant interval $I$, which contain the periodic point. For single letter cycles ($n=1$) these are the base intervals $I_{s}$, $s=0,\ldots,M$ in which the invariant interval is partitioned by the $M$ critical points of the map $f$. For any admissible cycle of topological length $n>1$ we:
\begin{enumerate}
 \item Generate all possible symbolic sequences of points on the cycle (as above).
 \item For cycle point with sequence $s_0s_1\ldots s_n$, compute $I_{s_0s_1\ldots s_n}$ recursively by application of the identity\rf{gilmore2003}
 \[
  I_{s_0s_1\ldots s_n}= I_{s_0}\cap f_{s_{0}}^{-1}(I_{s_1\ldots s_n})\,,
 \]
 where $f_{s_{0}}^{-1}$ is the inverse of the map $f$ restricted to the base interval $I_{s0}$ (in which $f$ is monotonic).
 \item Take the midpoint of $I_{s_0s_1\ldots s_n}$ as a guess for the cycle point $s_0s_1\ldots s_n$.
 \item Optionally, refine the guess by solving $f^n(q)=q$ through bisection or, for longer cycles, by setting up a multi-point shooting method~\cite{DasBuch}.
\end{enumerate}

Each admissible cycle of the return map can then be used to construct an initial ``guess'' solution, to be refined by an iterative method until it converges to the corresponding cycle of the continuous time flow (if the latter exists). 

\begin{enumerate}
 \item For each cycle point of the return map compute the corresponding full space point, through the the interpolating function $\mathcal{I}:\, q\mapsto x$, which maps LLE coordinates to ambient space coordinates $x\in R^d$ (Fourier modes $a_i$ or the case of KSe). Note that in the case of 1D KSe embedding we construct $d$ 1D interpolation functions $\mathcal{I}_i:\, q_1\mapsto a_i$, $i=1,\ldots d$. For the 2D embedding case of KSe with $\nu=0.02973$, we only needed to perform interpolation for points on $E_1$ and $E_2$ which are well parameterized by the first LLE coordinate $q_1$; this approach could be generalized to interpolate from $q$ to $a_i$, if needed.
 \item Use these approximate cycle points to set up a multiple-shooting or variational method to detect periodic orbits\rf{DasBuch}. In our implementation the collocation solver \texttt{solve\_bvp} build into \texttt{scipy}\rf{2020SciPy-NMeth} is used and solutions are detected in the full space, imposing the boundary condition $u(x,t+T_p)=u(x,t)$, with $T_p$ undetermined. In order to prepare an initial time-periodic guess the Poincar\'e section cycle points are used as initial conditions which are integrated until the next intersection with the Poincar\'e section. Special care is taken to lift points from the fundamental domain to the full space in order to ensure continuity of the initial guess. The guess is refined to a prescribed error tolerance, here $10^{-6}$.
 \item Once a periodic orbit is found we check that it corresponds to the expected label by computing its intersections with the Poincar\'e section, mapping all points to the fundamental domain and computing the transform of these points to LLE variables $q$ (this functionality is provided by the LLE implementation in \texttt{scikit-learn} without requiring recomputing the transform for all points in the dataset). 
\end{enumerate}

\section*{References}

%

\end{document}